\documentclass{conm-p-l}
\usepackage{amssymb}
\usepackage{amscd}
\newtheorem{theorem}{Theorem}[section]
\newtheorem{lemma}[theorem]{Lemma}

\theoremstyle{definition}

\newtheorem{example}[theorem]{Example}
\newtheorem{proposition}[theorem]{Proposition}
\newtheorem{cor}[theorem]{Corollary}

\theoremstyle{remark}
\newtheorem{remark}[theorem]{Remark}

 \newcommand{\Matrix}[4]{
  \left( \begin{array}{ccc}
   #1 & #2 \\
   #3 & #4 \\
  \end{array} \right) }
\newcommand{\uX}{{\underline X}}

\newcommand{\cX}{{\mathcal X}}
\newcommand{\cG}{{\tau_G}}
\newcommand{\cGr}{{\tau_{cG}}}
\newcommand{\cGg}{{\cGr^\bot}}

\newcommand{\cC}{{\mathcal C}}
\newcommand{\cI}{{\mathcal I}}

\newcommand{\mF}{{\mathbb F}}
\newcommand{\mT}{{\mathbb T}}
\newcommand{\cSi}{{R\mbox{-Sing}}}
\newcommand{\cSm}{{R\mbox{-Small}}}
\newcommand{\ZZ}{{\mathbb Z}}
\newcommand{\QQ}{{\mathbb Q}}
\newcommand{\NN}{{\mathbb N}}
\newcommand{\ess}{\trianglelefteq}
\newcommand{\Rad}[1]{{\mathrm{Rad}}\:(#1)}
\newcommand{\Soc}[1]{{\mathrm{Soc}}\:(#1)}

\newcommand{\Hom}[2]{{\mathrm{Hom}}\:(#1,#2)}
\newcommand{\Jac}[1]{{\mathrm{Jac}}\:(#1)}

\renewcommand{\Im}[1]{{\mathrm{Im}}\:(#1)}

\numberwithin{equation}{section}

\begin{document}

\title{On the Splitting of the Dual Goldie Torsion Theory}

\author{Christian Lomp}
\address{ Departamento de Matem\'atica Pura, Universidade do Porto,
4050 Porto, Portugal
}
\email{clomp@fc.up.pt}
\thanks{This note was done with the support of "DAAD
Doktorandenstipendium HSP III" while the author was visiting
the Department of Mathematical Science at UW-Milwaukee 1997/98.
He would like
to thank the department for their support and hospitality. Moreover
he would like to express his thanks
 to Ay\c{s}e \c{C}. \"Ozcan for supplying him with her paper and her
 thesis and
 to Mark L. Teply for his valuable comment and for Proposition \ref{teply}.}



\date{today}


\begin{abstract} The splitting of the Goldie (or singular) torsion theory
has been extensively studied. Here we determine an appropriate
dual Goldie torsion theory, discuss its splitting and answer in
the negative a question proposed by \"Ozcan and Harmanc{$\imath$}
as to whether the splitting of the dual Goldie torsion theory
implies the ring to be quasi-Frobenius.
\end{abstract}
\maketitle

\section{Introduction}
Let $R$ be an associative ring with unit and $M$ a left
$R$-module. A submodule $N$ of $M$ is called {\it essential in
$M$} (denoted by $N \ess M$) if $N\cap L \neq 0$ holds for every
non-zero submodule $L$ of $M$. A left
$R$-module $M$ is called {\it singular} if there are left
$R$-modules $L$ and $K$ such that $K$ is an essential submodule of
$L$ and $M \simeq L/K$. We denote the class of singular left
$R$-modules by $\cSi$ and the sum of all singular submodules of a
module $M$ by $Z({_RM}) := trace({\cSi}, M)$, where the trace of a
class of modules into a module is defined by $trace(\cX, M) :=
\sum \{\Im{f} : f \in \Hom{X}{M}, X \in \cX \}$.
 The lattice-theoretical dual of the notion of essential submodules
are {\it small }  submodules $N$ of $M$ (denoted by $N \ll M$)
that have the property that
 $N + L \neq M$ holds for every proper submodule $L$ of $M$. Dual to a
 singular module, a
left $R$-module $M$ is called {\it small} if there are left
$R$-modules $K$ and $L$ such that $K$ is a small submodule of $L$
and $M \simeq K$. A definition of small and singular objects of abelian
categories was given by B.Pareigis (see \cite{pareigis}).
It is easily verified that $M$ is a small left
$R$-module if and only if it is a small submodule of its injective
hull $E(M)$ (see for example \cite[Theorem 1]{leonard}).
We denote the class of small left
$R$-modules by $\cSm$ and denote the sum of all small
submodules of a module $M$ by $Z^*({_RM}):=trace({\cSm},M)$. The
class of small modules is closed under submodules, homomorphic
images and finite direct sums.
Note that a simple left $R$-module is either singular or projective and
either small or injective.
Small and singular modules occur in decomposition theorems of modules
over quasi Frobenius rings (QF-rings).
\begin{theorem}[Rayar {\cite[Theorem 7, Corollary 11]{Rayar1}}, Oshiro \cite{oshiro}]
Let $R$ be a ring. The following are equivalent:
\begin{enumerate}
\item[(a)] $R$ is a QF-ring;
\item[(b)] every left $R$-module is a direct sum of a projective and a
small module;
\item[(c)] every left $R$-module is a direct sum of an injective and a singular
module;
\item[(d)] $\Soc{_RR} \subseteq \Soc{R_R}$ holds and
    \begin{enumerate}
        \item[(i)]
        every left $R$-module is a direct sum of an injective and a
        small module or
        \item[(ii)]
        every right $R$-module is a direct sum of a projective and a
        singular module.
    \end{enumerate}
\end{enumerate}
\end{theorem}

Recently \"Ozcan and Harmanc{$\imath$} proved in \cite{ozcan} that
modules over QF-rings have a direct sum decomposition into two
"orthogonal" classes of modules. In Section 2 we will realize
those two module classes as a torsion and torsion free class of a
hereditary torsion theory - the dual Goldie torsion theory - and
state some basic properties. Those rings whose dual Goldie torsion
theory is trivial or improper will be studied in Section 3.
Especially those rings $R$ that are small $R$-submodules of their
injective hull will be discussed. In Section 4 we will determine
when the dual Goldie torsion theory splits. In particular we will
characterize semilocal rings whose dual Goldie torsion theory is
splitting as those that cogenerate all injective simple left
$R$-modules. Hence semilocal left Kasch rings have this property.
We will give a list of classes of rings whose dual Goldie torsion
theory is splitting, but that are far from being QF. This answers
in the negative a question of \"Ozcan and Harmanc{$\imath$} as to
whether the splitting of the dual Goldie torsion theory implies
the ring to be QF (see \cite[pp 325]{ozcan}). For all unexplained
ring- and module-theoretical notations we refer to
\cite{wisbauer}; for all unexplained torsion-theoretical notations
we refer to \cite{golan}.

\section{Dual Goldie Torsion Theory}
Let $\cX$ be a class of left $R$-modules closed under ismorphisms and submodules.
Define the following classes:
\begin{equation*}
\begin{split}
 \mF(\cX) :=& \{ M\in R\mbox{-Mod  } | \forall X\in \cX : \Hom{X}{M}=0\}\\
 \mT(\cX) :=& \{ N\in R\mbox{-Mod  } | \forall M \in \mF(\cX) :
 \Hom{N}{M}=0 \}
\end{split}
\end{equation*}
Then $( \mT(\cX), \mF(\cX) )$ is a hereditary torsion theory (see
\cite[II1.3]{bican}), the smallest hereditary torsion theory such that all
$\cX$-modules are  torsion.
The torsion radical is $\tau_{\cX}(M) = trace( \mT(\cX), M)$ and
$trace(\cX, M) \ess \tau_{\cX}(M)$ holds. This kind of torsion theory has
also been extensively studied by Harmanc{$\imath$} and Smith, where
the corresponding torsion radical was denoted by ${\underline H}^*$ (see in \cite{harmanci}).

The following are easily checked:
\begin{equation} \label{lemma1}
\begin{split}
\mT(\cX) =& \{ N\in R\mbox{-Mod } | \forall U \subset V \subseteq N \mbox{ : } trace(\cX, V/U) \neq 0 \}\\
\mF(\cX) =&\{ M \in  R\mbox{-Mod  } | trace(\cX,M)=0 \}
\end{split}
\end{equation}

The {\it Goldie torsion theory} is defined as:
$ \cG = ( \mT({\cSi}), \mF({\cSi}) ).$
V.S. Ramamurthi defined a {\it dual Goldie torsion theory} as
$ \cGr = ( \mT({\cSm}), \mF({\cSm}) )$
and studied some of its properties in \cite{ramamurthi}.
We have $Z^*(M)=M\cap \Rad{E(M)} \ess \cGr(M)$ for every left $R$-module $M$.

The class of singular modules is closed under
arbitrary direct sums, but does not have to be closed under extensions
in general. In case of a nonsigular ring $\cSi = \mT({\cSi})$. The class
of small modules is, in general, neither closed under arbitrary direct sums
nor under extensions. A necessary and sufficient condition for $\cSm$ to
be closed under direct sums is that every injective left $R$-module has
a small radical (see \cite[Lemma 9]{Rayar1}).

Related to $\cGr$ are two classes of modules introduced by
\"Ozcan and Harmanc$\imath$ in \cite{ozcan}:
\begin{equation*}
\begin{split}
\uX   :=& \{ M\in R\mbox{-Mod  } | Z^*(M)=0  \} \\
\uX^* :=& \{ M\in R\mbox{-Mod  } | \forall U\subset V\subseteq M :
Z^*(V/U) \neq 0\}
\end{split}
\end{equation*}
From (\ref{lemma1}) we see that $\cGr = (\uX^*,\uX)$. \"Ozcan and
Harmanc{$\imath$} showed that over a QF-ring every left $R$-module
is a direct sum of an $\uX$-module and an ${\uX}^*$-module. In
other words they showed that $\cGr$ splits over QF-rings. They
raised the question, ``Is a ring $R$ whose dual Goldie torsion
theory $\cGr$ is splitting a QF ring ?´´ (see \cite[pp
325]{ozcan}). It is not difficult to see that this question has a
negative answer, and we will discuss the splitting of $\cGr$ in
the sequel.

A different approach to defining a dual Goldie torsion theory was
proposed by A.I. Generalov in \cite{generalov}. Let $\cX$ be a
class of left $R$-modules closed under ismorphisms and submodules.
Define the following classes:
\begin{equation*}
\begin{split}
{\mT^\bot}(\cX) :=& \{ M\in R\mbox{-Mod  } | \forall X\in \cX : \Hom{M}{X}=0 \}\\
{\mF^\bot}(\cX) :=& \{ N\in R\mbox{-Mod  } | \forall M \in {\mT^\bot}(\cX) : \Hom{M}{N}=0 \}
\end{split}
\end{equation*}
Then $({\mT^\bot}(\cX), {\mF^\bot}(\cX))$ is a torsion theory, not necessarily hereditary.
Denote the reject of $\cX$ in $M$ by
$reject(M,\cX) := \bigcap \{\ker{g} : g \in \Hom{M}{X}, X \in \cX \}$.
Then we see as above:
\begin{equation}\label{lemma2}
\begin{split}
{\mT^\bot}(\cX) =& \{ M \in  R\mbox{-Mod  } | reject(M,\cX)=M \}\\
{\mF^\bot}(\cX) =& \{ N \in  R\mbox{-Mod  } | \forall 0\subset U \subseteq M
\mbox{ : } reject(U, \cX) \neq U \}
\end{split}
\end{equation}
The torsion class ${\mT^\bot}(\cSm)$ has the property that $M\in
{\mT^\bot}(\cSm)$ if and only if every factor module of $M$ is
$\cGr$-torsion free. We denote this torsion theory by $\cGg$, for
if $\cGg$ is hereditary, then it is the pseudo-complement of
$\cGr$ in the lattice of hereditary torsion theories $R$-tors (see
\cite[pp. 280]{golan}).

A characterization of $\cGg$-torsion modules was given in
\cite{generalov} and \cite{lomp}:

\begin{proposition}[Lomp {\cite[5.1]{lomp}}, Generalov {\cite[Proposition 3]{generalov}}]

Let $R$ be a ring and let $M$ be a left $R$-module. Then the
following conditions are equivalent:
\begin{enumerate}
\item[(a)] $M$ is $\cGg$-torsion;
\item[(b)] $\forall K \in R$-Mod and $0 \neq f: M \rightarrow K$, $\Im{f}
\not\ll K$;
\item[(c)] $\forall K \in R$-Mod and $0 \neq f: M \rightarrow K$, $\Im{f}$ is
coclosed in $K$, i.e. $\Im{f}/L \not\ll K/L$ for all $L\subset \Im{f}$ ;
\item[(d)] for every factor module $L$ of  $M$ and diagram in $R$-Mod
$$\begin{CD}
0 @>>> A @>i>> B \\
@. @VfVV \\
@. L
\end{CD}$$
there exists a factor module $\pi: L \rightarrow \tilde L$
and a $g\in \Hom{B}{\tilde L}$ such that $ig=f\pi$ holds.
\end{enumerate}
\end{proposition}

The following result is easily verified:
\begin{lemma}
$\mF(\cSm) = {\mT^\bot}(\cSm)$ if and only if
$\cGr$ is cohereditary and $\cGg$ is hereditary.
\end{lemma}

The torsion radical of $\cGg$ can also be defined via an inductive
construction of a preradical $\rho$ given by Generalov: for any
module $M$ let $\sigma^1(M): =reject(M,\cSm)$, for any ordinal
number $\alpha$ define $\displaystyle {\sigma^{\alpha+1}(M) : =
reject(\sigma^\alpha(M),\cSm)}$ and if $\alpha$ is a limit
ordinal, let $\sigma^{\alpha}(M) : = \bigcap_{\beta < \alpha}
\sigma^{\beta}(M)$. Then there exists an ordinal number $\gamma$
such that $\sigma^\gamma(M)=\sigma^{\gamma +1}(M)$. Hence define a
preradical by $\rho(M):=\sigma^{\gamma}(M)$. Generalov showed that
$\rho$ is an idempotent preradical, that defines the torsion class
${\mT^\bot}(\cSm)$, and called it the dual Goldie torsion theory.

\section{V-Rings, Small Rings and Almost Small Rings}
All torsion theories are considered as torsion theories of left
$R$-modules unless otherwise specified.
We will use the torsion theoretical notion for $\xi$ to denote the trivial
torsion theory (all non-zero modules are torsion free ) and $\chi$ to denote
the improper torsion theory (all modules are torsion).

\subsection{V-Rings}
We first examine when $\cGr$ becomes trivial.
\begin{proposition}
Let $R$ be a ring.
Then $R$ is a left $V$-ring if and only if $\cGr=\xi$
if and only if $_RR$ is $\cGg$-torsion.
\end{proposition}

\begin{proof}
Simple modules are either small or injective (see \cite{Rayar1}).
If $\cGr=\xi$, then all simple left $R$-modules are injective and
hence $R$ is a left $V$-ring. On the other hand, if $R$ is a left
$V$-ring, then $\Rad{M}=0$ for all left $R$-modules $M$. Hence
there are no non-zero small left $R$-modules and $\cGr=\xi$.
Obviously $\cGr=\xi \Leftrightarrow \cGg=\chi \Leftrightarrow
{_RR}$ is $\cGg$-torsion.
\end{proof}

If $\cGr=\xi$, then $\cGr$ is trivially splitting; e.g. for a
direct product of fields $\cGr$ is splitting. The following
example of a $\cGr$-torsion free ring $R$ with $\cGr \neq \xi$
shows that $R$ being $\cGr$-torsion free is not sufficient for $R$
to be a $V$-ring.

\begin{example}\label{example_inf_vs} The endomorphism ring $S$ of an infinite dimensional
 vector space $V_k$ over a field $k$ is a von Neumann regular
left self-injective ring, but not a left $V$-ring (since $_SV$ is
a simple, non-injective $S$-module and hence small, see
\cite[23.6]{wisbauer}). As $S$ is left self-injective, we have
$\Jac{S} = Z^*(_SS)$. Moreover $\Jac{S}=0$, as $S$ is von Neumann
regular. Thus $S$ is $\cGr$-torsion free, but $\cGr \neq \xi$.
\end{example}

\subsection{Small Rings}
We will now examine when $\cGr$ becomes improper. A ring $R$ is
called {\it left small} if $_RR$ is a small module; e.g. $\ZZ$ is
a small ring as it is small in $_\ZZ\QQ$. A left small ring $R$ is
$\cGr$-torsion as left $R$-module. Hence every left $R$-module is
$\cGr$-torsion and we have $\cGr = \chi$.

\begin{proposition}[Ramamurthi {\cite[3.3]{ramamurthi}}, Pareigis {\cite[Satz 4.8]{pareigis}}]
Let $R$ be a ring and let $E(R)$ be the injective hull of $_RR$.
Then the following conditions are equivalent:
\begin{enumerate}
\item[(a)]$R$ is a left small ring;
\item[(b)]$\Rad{M}=M$ for every injective left
$R$-module $M$;
\item[(c)] $\Rad{E(R)}=E(R)$.
\end{enumerate}
\end{proposition}

Harada showed that a commutative ring is small in its classical quotient ring
if and only if every maximal ideal contains a regular element (\cite[Theorem 2]{harada}).

\begin{theorem} \label{primitive}
Let $R$ be a ring such that every left primitve ideal contains
a regular element. Then $R$ is left small.
\end{theorem}

\begin{proof}
Assume $R$ is not small in its injective hull $E:=E(R)$. Then
there exists a submodule $N\subseteq E$ such that $R+N=E$ as left
$R$-modules. Hence $E/N$ is cyclic and has a maximal submodule
$M/N$. The left annihilator $I:=l(E/M)$ of the left simple
$R$-module $E/M$ contains a regular element
 $x$ by hypothesis. As $E$ is a divisible left $R$-module, we have
$x(E/M)=0 \Rightarrow E=xE\subseteq M$, which is a contradiction.
Hence $R$ has to be small in $E(R)$.
\end{proof}

As a corollary we get the following known results:

\begin{cor}\label{cor_prim}
Proper integral domains (Pareigis, \cite[Folgerung 5.3]{pareigis}),
rings whose Jacobson radical contains a regular element
(Ramamurthi, \cite[Proposition 3.4]{ramamurthi}) and
prime left Goldie rings that are not left primitive
(\"Ozcan, \cite[2.2.5]{cigdemThesis}) are small rings.
\end{cor}

Recall that a left $R$-module $M$ is torsionfree (in the usual
sense) if $rm\neq 0$ holds for every $m\in M$ and $r\in R$. The
following result shows that small rings bring some restrictions.

\begin{proposition}
Let $R$ be a left small ring, $E$ an injective left $R$-module and
$\lambda$ the cardinality of an generating set of $E$. Then
$\lambda \geq \aleph_0$ and $R^{(\lambda)}$ is not a small left
$R$-module. Moreover the cardinality of any independent family of
submodules of a finitely generated torsionfree left $R$-module $M$
is bounded by $\lambda$. If $\lambda = \aleph_0$, every finitely
generated torsionfree left $R$-module has finite Goldie dimension.
\end{proposition}

\begin{proof}
First note that $E$ is generated by $R^{(\lambda)}$ and that an
injective left $R$-module cannot be small. Therefore $\lambda$ has
to be infinite and $R^{(\lambda)}$ cannot be a small module. If we
have an independent family of cardinality $\lambda'$ of cyclic
submodules of a finitely generated torsionfree left $R$-module
$M$, then $R^{(\lambda')}$ embeds into $M$. If $\lambda' \geq
\lambda$, then $R^{(\lambda')}$ is not small. Hence $M$ is not
small, which contradicts the fact that $M$ is a small module as it
is a homomorphic image of a finite direct sum of copies of the
small module $_RR$. Thus $\lambda' < \lambda$.
\end{proof}

\begin{remark}
Small rings form a class of rings where small modules are not
closed under arbitrary direct sums. Contrary to this class are
left max rings. A ring $R$ is called a left {\it max ring} if
every left $R$-module has a maximal submodule, or equivalently, if
every left $R$-module has a small radical. By M.Rayar's result
mentioned above, $\cSm$ is closed under arbitrary direct sums over
a max ring $R$. It is well-known that $R$ is left perfect if and
only if $R$ is a semilocal left max ring. Hence we already
encountered a huge class of rings, namely the small rings,
 whose dual Goldie torsion theory (trivially) splits, but that cannot be
 perfect (nor QF) answering in the negative a question of \"Ozcan and Harmanc{$\imath$}
as to whether the splitting of $\cGr$ implies the ring to be QF.
 As an additional observation together with \ref{primitive}, we have that
 a left max domain has to be left primitive.
\end{remark}

\subsection{Almost Small Rings}
Let us call a ring $R$ {\it left almost small} if $\cGr=\chi$.
Obviously every left small ring is left almost small, and a ring
$R$ is left almost small if and only if $_RR$ is $\cGr$-torsion.
So $\cGr$ trivially splits for left almost small rings.

\begin{proposition}\label{proplocal}
Any local ring is either a division ring or else left and right
almost small.
\end{proposition}

\begin{proof}
Let $R$ be a local ring. If $R/\Jac{R}$ is injective, then $R$ is
a left $V$-ring. Hence $\Jac{R}=0$, which implies $R$ is a
division ring. If $R$ is not a division ring, then $R/\Jac{R}$
must be a small left (or right) $R$-module. Thus $R$ is
$\cGr$-torsion, as $\Jac{R}$ is small and $\cGr$-torsion modules
are closed under extensions. Hence $R$ is left and right almost
small.
\end{proof}

\section{The splitting of the dual Goldie torsion theory.}

There are many example of rings where the dual Goldie torsion
theory splits. M.Rayar proved in \cite[Proposition 1]{Rayar2} that
a direct product of a family of proper integral domains is a small
ring. Hence any product of integral domains is a direct product of
a small and a $V$-ring; thus $\cGr$ splits. Take, for instance,
$R=\ZZ^\NN \times \QQ^\NN$.

We will need the next technical lemma to investigate torsion theories
of finite products of rings.

\begin{lemma} Let ${R_1}, \ldots ,{R_n}$ be rings and denote by $R:={R_1}
\times\cdots \times {R_n}$ their direct product with componentwise
multiplication. For all $i$, let $\cX_i \subseteq R_i-Mod$ be
classes of left $R_i$-modules closed under submodules and
isomorphisms and set $\cX := \bigoplus_{i=1}^n \cX_i$. Then
$\mF(\cX) = \bigoplus_{i=1}^n \mF(\cX_i)$ and $\mT(\cX) =
\bigoplus_{i=1}^n \mT(\cX_i)$ (where the $\mF(\cX_i)$, resp.
$\mT(\cX_i)$, are formed in $R_i-Mod$).
\end{lemma}

\begin{proof}
Note that the unit of $R$ is $1 = (1_{R_1}, \ldots, 1_{R_n})$ and
that $R=\bigoplus_{i=1}^n R_i$ as a left $R$-module. Hence for all
$M \in R-Mod$ we have $M = \sum_{i=1}^n M_i$ with $M_i := R_iM$.
Every left $R_i$-module $M$ becomes a left $R$-module by $(r_1,
\ldots, r_n) \cdot m := r_im$. Note that $R_i \cdot M_j = 0$ if
$i\neq j$. Therefore every $M \in R-Mod$ can be decomposed as $M =
\bigoplus_{i=1}^n M_i$. Let $M \in R_i-Mod$ and $N \in R-Mod$ such
that $N_i=0$; then $Hom_R(N,M)=0$.
Now let $N, M \in R-Mod$: $$Hom(M,N) \simeq
\bigoplus_{i=1}^n\bigoplus_{j=1}^n Hom(M_i, N_j) =
\bigoplus_{i=1}^n Hom(M_i, N_i).$$ Hence $trace(\cX, M)
=\bigoplus_{i=1}^n trace(\cX_i, M_i)$ holds, which implies
$\mF(\cX) = \bigoplus_{i=1}^n \mF(\cX_i)$ and $\mT(\cX) =
\bigoplus_{i=1}^n \mT(\cX_i)$.
 \end{proof}

Thus as a direct
consequence we get the following corollary (note that
${\cSm} = \bigoplus_{i=1}^n {R_i\mbox{-Small}}$ holds).
Let us denote by $_{R_i}\tau$ a torsion
theory in $R_i$-tors.

\begin{cor} \label{cor_prod}
Let $R=R_1 \times \cdots \times R_n$. Then ${_R\cGr} = {\left(
_{R_1}\cGr \right)} \times \cdots \times {\left( _{R_n}\cGr
\right)}$.
\end{cor}

As a consequence of Corollary \ref{cor_prod} and the facts that
every commutative semiperfect ring is  a direct product of local
rings and every local ring is either a division ring or almost
small by Proposition \ref{proplocal}, we observe:

\begin{cor}
Any commutative semiperfect ring is a direct product of an almost
small ring and a semisimple ring.
\end{cor}

A torsion theory $\tau$ is called {\it cohereditary} if the class
of $\tau$-torsionfree modules is closed under homomorphic images.
Example \ref{example_inf_vs} shows that $\cGr$ does not have to be
cohereditary. With the help of the above Corollary \ref{cor_prod}
we can determine when $\cGr$ splits for rings with cohereditary
dual Goldie torsion theory:

\begin{theorem}\label{thm1} Let $R$ be a ring. Then the following
conditions are equivalent:
\begin{enumerate}
 \item[(a)] $\cGr$ is splitting and cohereditary;
 \item[(b)] $R \simeq T \times S$ with $T$ a left almost small and $S$ a
left V-ring;
\end{enumerate}
\end{theorem}

\begin{proof}
$(a) \Rightarrow (b)$ Let ${_RR} = {_RT} \oplus {_RS}$ with $_RT =
\cGr(_RR)$  and $_RS$ being $\cGr$-torsion free. We have
$\Hom{T}{S}=0$, and thus $TS=0$ implies that $T$ is an ideal in
$R$. As the class of $\cGr$-torsion free modules is closed under
homorphic images, we have $\Hom{S}{T}=0$ and thus $ST=0$. Hence
$T$ and $S$ are ideals with $T \cap S = 0$. Thus $R = T \times S$,
and as $\cGr$ is cohereditary, all $\cGr$-torsion left $R$-modules
are generated by $_RT$ and all $\cGr$-torsion free left
$R$-modules are generated by $_RS$. Thus $T$ is left almost small
and $S$ a left V-ring.\\ $(b) \Rightarrow (a)$ By Corollary
\ref{cor_prod} we have $(\mT(\cSm),\mF(\cSm))= (T-Mod, S-Mod)$,
and every left $R$-module $M$ can be decomposed into $M=TM \oplus
SM$. Thus $\cGr$ is splitting. $\cGr$ is cohereditary since
$\mF(\cSm) = S\mbox{-Mod}$.
\end{proof}

 A module $M$ is called {\it cosemisimple}
if $\Rad{M/K}=0$ for all $K \subset M$. All $\cGr$-torsion free
modules are cosemisimple in case that $\cGr$ is cohereditary. In
the next proposition we will characterize when $\cGr$ is
cohereditary. Recall that a module $M$ is called $\tau$-{\it
injective} with respect to a hereditary torsion theory $\tau$ if
$M$ is injective with respect to all short exact sequences $0
\rightarrow L \rightarrow N \rightarrow K \rightarrow 0$ such that
$K$ is $\tau$-torsion.

\begin{proposition} Let $R$ be a ring. Then the following conditions are equivalent:
\begin{enumerate}
 \item[(a)] $\cGr$ is cohereditary;
 \item[(b)] all $\cGr$-torsion free left $R$-modules are $\cGr$-injective;
 \item[(c)] $E(M)/M$ is $\cGr$-torsion free for all $\cGr$-torsion free left
$R$-modules $M$.
\end{enumerate}
If one of the above conditions hold, then all $\cGr$-torsion free
left $R$-modules are cosemisimple.
\end{proposition}

\begin{proof}
$(a) \Rightarrow (c)$ is trivial; $(c) \Leftrightarrow (b)$ holds
for any torsion theory (see \cite{wisbauer96}). $(b) \Rightarrow
(a)$ Since $\cGr$ is hereditary, every submodule $N$ of a
$\cGr$-torsion free module $M$ with $M/N$ being $\cGr$-torsion
satisfies $M \simeq N \oplus M/N$ and so $N=M$.
 \end{proof}

For semilocal rings $R$ the dual Goldie torsion theory becomes
rather simple. Let $R$-Simp, and let $\cC$, resp. $\cI$, denote
the class of simple, semisimple, resp. injective, left
$R$-modules. Before we state a more general result about the
splitting of $\cGr$. We note that the following equivalences hold:
$$\mF{(\cSm)} \subseteq \cC \Leftrightarrow \mF{(\cSm)} \subseteq
\cI \Leftrightarrow \mF{(\cSm)} = \cC \cap \cI.$$ Moreover in this
situation $\cGr$ is cohereditary and the dual Goldie torsion
theory is nothing but the torsion theory cogenerated by the
injective simple left $R$-modules; i.e. $\cGr =
\chi\{R\mbox{-Simp} \cap \cI \}$. We have this situation when
$R/Z^*(_RR)$ is semisimple; e.g., $R$ is semilocal.
\begin{theorem} \label{thm2}
Let $R$ be a ring with $R/Z^*(_RR)$ semisimple. Then the following
conditions are equivalent:
\begin{enumerate}
\item[(a)] $\cGr$ is splitting;
\item[(b)] $R \simeq T \times S$ with $T$ a left almost small and $S$ a semisimple ring;
\item[(c)] $\cGr$ is stable ( i.e. the $\cGr$-torsion class is closed under essential
extensions);
\item[(d)] all $\cGr$-torsion free modules are projective;
\item[(e)] all injective simple left $R$-modules are projective;
\item[(f)] $\cG \leq \cGr$ (i.e, $\cGr$ is a generalization of $\cG$);
\item[(g)] all singular left $R$-modules are $\cGr$-torsion;
\item[(h)] $R$ cogenerates all injective simple left $R$-modules.
\end{enumerate}
\end{theorem}

\begin{proof}
If $R/Z^*(_RR)$ is semisimple, then all $\cGr$-torsion free
modules are semi-simple as $Z^*(_RR)M=0$. By the above we have
$\mF{(\cSm}) = \cC \cap \cI$; thus $\cGr$ is cohereditary. Hence
Theorem \ref{thm1} gives us $(a) \Leftrightarrow (b)$. Note that
$(b) \Leftrightarrow (c)$ by \cite[Prop. 5.10]{golan},
$(c) \Leftrightarrow (d)$ by \cite[E5.10]{golan}, and $(d)
\Leftrightarrow (e)$ is obvious. To show $(d) \Rightarrow (f)$ let
$M$ be $\cGr$-torsion free; then $Z(M)$ being projective implies
$Z(M)=0$ and $M$ is non-singular. Hence all $\cGr$-torsion free
modules are $\cG$-torsion free: $\cG \leq \cGr$. The implication
$(f) \Rightarrow (g)$ is trivial. To prove $(g) \Rightarrow (e)$,
note that an injective simple module cannot be small; therefore it
cannot be singular by hypothesis and must be projective. Finally,
$(e) \Leftrightarrow (h)$ is obvious.
 \end{proof}

\begin{remark}
Semilocal rings satisfy the hypothesis of Theorem \ref{thm2}. A
ring $R$ is called a {\it left Kasch} ring if it cogenerates all
left simple $R$-modules. Of course any left cogenerator ring
(e.g., left PF-ring or QF-ring) is a left Kasch ring, but the
converse does not hold in general. (A ring is called left PF-ring
if it is left self-injective and left Kasch.)
\end{remark}

 As a corollary of Theorem \ref{thm2} we note:

\begin{cor} For any semilocal left Kasch ring, $\cGr$ splits.
\end{cor}

\begin{remark}
A result by Faith and Menal \cite{FaithMenal} states that a left
Kasch ring with finite right uniform dimension is semilocal.
Another result by G\'omez Pardo and Yousif \cite{GomezPardoYousif}
states  that a left Kasch, right CS ring is semiperfect. Hence for
those rings $\cGr$ splits.
\end{remark}

A sufficient and necessary condition for a semilocal ring to be
left Kasch is the following result that was pointed out to the
author by Mark Teply. The left (resp. right) annihilator of an
ideal $I$ is denoted by $l(I)$ (resp. $r(I)$).

\begin{proposition} \label{teply}
A semilocal ring $R$ with $l(r(\Jac{R}))=\Jac{R}$ is a left Kasch ring.
\end{proposition}

\begin{proof}
Assume that $l(r(\Jac{R}))=\Jac{R}$ holds and that $_RE$ be a
simple module with $E\simeq (Ra +\Jac{R})/\Jac{R} \simeq
Ra/(Ra\cap \Jac{R}) $ for an $a\in R$. Since $a\notin \Jac{R}$,
there exists by hypothesis $b\in r(\Jac{R})=\Soc{_RR}$ with $ab
\neq 0$. Thus we have an isomorphism $\phi : E \rightarrow Rab$
with $\phi(r\bar{a}):=rab$. Hence $R$ cogenerates all simple left
$R$-modules; $R$ is a left Kasch ring. The converse is true for
every left Kasch ring.
\end{proof}

\begin{remark}
Summarizing we have that $\cGr$ splits for the following clases of
rings:
\begin{itemize}
\item semilocal left Kasch rings;
\item left Kasch rings with finite right Goldie dimension;
\item left Kasch, right CS-rings;
\item left V-rings;
\item local rings;
\item prime left Goldie rings, not left primitive;
\item products of integral domains;
\item commutative semiperfect rings;
\item commutative noetherian semilocal rings;
\item any finite direct product of the rings above.
\end{itemize}
Note that commutative noetherian rings are ``stable´´ rings; i.e.
every hereditary torsion theory is stable. By Theorem \ref{thm2}
we know that $\cGr$ splits for semilocal stable rings.
\end{remark}
\begin{remark}
There are (local) rings with $\cGr$ splitting and additional
reasonable good properties (like 'injective cogenerator' or
'artinian') that are still not QF. For instance  F.Dischinger and
W.M\"uller give a famous example of a local left PF-ring that is
not right PF (see \cite{dischinger}), and Bj\"orks gives an
example of a local two-sided artinian ring with $\Jac{R}^2=0$ that
is not left self-injective (see \cite{bjoerk}).

On the other hand, even a two-sided artinian condition need not
imply the splitting of $\cGr$. Let $k$ be a field and
$R:=\Matrix{k}{k}{0}{k}$. Then $R$ is two-sided artinian. The left
simple $R$-module $R/{\Soc{_RR}}$ is injective but not projective,
and hence by Theorem \ref{thm2} $\cGr$ does not split.
\end{remark}

\begin{proposition} Let $R$ be a left artinian ring with
${\Soc{R_R}} \subseteq {\Soc{_RR}}$. Then $\cGr$ splits.
\end{proposition}

\begin{proof} Any left singular $R$-module $M$ is a
$R/{\Soc{_RR}}$-module and so by hypothesis is also a
$R/{\Soc{R_R}}$-module. By \cite[Theorem 3]{Rayar1} $M$ is a left
small $R$-module; by Theorem \ref{thm2}(h) $\cGr$ splits.
 \end{proof}

We would like to end this note with a remark:
In view of Theorem \ref{thm1}, is there  an example of a ring such
that $\cGr$ is splitting but not cohereditary ?

\bibliographystyle{amsalpha}

\end{document}